\input amstex
\input Amstex-document.sty

\pageno 163

\def\shorttitle{Tamagawa Number Conjecture for zeta Values}
\def\shortauthor{Kazuya Kato}

\topmatter
\title\nofrills{\boldHuge Tamagawa Number Conjecture for zeta Values}
\endtitle

\author \Large Kazuya Kato* \endauthor

\thanks *Department of Mathematical Sciences, University of Tokyo,
Komaba 3-8-1, Meguro, Tokyo, Japan. E-mail: kkato\@ms.u-tokyo.ac.jp
\endthanks

\abstract\nofrills \centerline{\boldnormal Abstract}

\vskip 4.5mm

{\ninepoint Spencer Bloch and the author formulated a general conjecture (Tamagawa number conjecture) on the
relation between values of zeta functions of motives and arithmetic groups associated to motives. We discuss this
conjecture, and describe some application of the philosophy of the conjecture to the study of elliptic curves.

\vskip 4.5mm

\noindent {\bf 2000 Mathematics Subject Classification:} 11G40.

\noindent {\bf Keywords and Phrases:} zeta function, Etale cohomology, Birch Swinnerton-Dyer conjecture.}
\endabstract
\endtopmatter

\def\bC{\bold C}
\def\bR{\bold R}
\def\bQ{\bold Q}
\def\bZ{\bold Z}
\def\bF{\bold F}
\def\cF{\Cal F}
\def\L{\Lambda}

\def\dt{\text{det}}

\def\Gal{\text{Gal}}
\def\Spec{\text{Spec}}
\def\RG{R\Gamma}
\def\1/p{{\bold Z}[\frac{1}{p}]}
\font\boldmath = cmmib10 scaled 1000

\document

\baselineskip 4.5mm \parindent 8mm

Mysterious relations between zeta functions and various arithmetic groups have been important subjects in number
theory.

\noindent (0.0) \ \ \ \ \ \ \ \ zeta functions $\leftrightarrow$ arithmetic groups.

A classical result on such relation is the  class number formula
discovered in 19th century, which relates zeta functions of number
field to ideal class groups and unit groups. As indicated in
(0.1)--(0.3) below, the formula of Grothendieck expressing the
zeta functions of varieties over finite fields by etale cohomology
groups, Iwasawa main conjecture proved by Mazur-Wiles, and Birch
and Swinnerton-Dyer conjectures for abelian varieties over number
fields, considered in 20th century, also have the form (0.0).

\noindent (0.1) \  Formula of Grothendieck.

\ \ \ \ \ \ \ \ zeta functions $\leftrightarrow$ etale cohomology groups.

\noindent (0.2) \ Iwasawa main conjecture.

\ \ \  \ \ \ \ \ zeta functions, zeta elements $\leftrightarrow$ ideal class groups, unit groups.

\noindent (0.3) \ Birch Swinnerton-Dyer conjectures (see 4).

\ \ \ \ \ \ \ \ zeta functions $\leftrightarrow$ groups of rational points, Tate-Shafarevich groups.

Here in (0.2), ``zeta elements" mean cyclotomic units which are
units in cyclotomic fields and closely related to zeta functions.
Roughly speaking, the relations (often conjectural) say that the
order of zero or pole of the zeta function at an integer point is
equal to the rank of the related finitely generated arithmetic
abelian group (Tate, the conjecture (0.3), Beilinson, Bloch, ...)
and the value of the zeta function at an integer point is related
to the order of the related arithmetic finite group.

In [BK], Bloch and the author formulated a general conjecture
on (0.0) (Tamagawa number conjecture for motives).
Further generalizations of Tamagawa number conjecture by
Fontaine, Perrin-Riou, and the author [FP],
[Pe$_1$]  [Ka$_1$], [Ka$_2$] have the
form

\noindent (0.4) \ \ \ \ \ zeta functions (= Euler products, analytic)

\ \ \ \ \ $\leftrightarrow$ zeta elements (= Euler systems,
arithmetic)

\ \ \ \ \ $\leftrightarrow$   arithmetic groups.

\noindent Here the first $\leftrightarrow$  means that zeta
functions enter the arithmetic world transforming themselves into
zeta elements, and the second  $\leftrightarrow$ means that zeta
elements generate ``determinants" of certain etale cohomology
groups.

 The aim of this paper is to discuss (0.4) in
an expository style. We review (0.1) in \S1, and then in
\S2, we
describe the generalized Tamagawa number conjecture (0.4),
the relation with (0.2), and an application of the
philosophy (0.4) to (0.3).

In this paper, we fix a prime number $p$. For a commutative
ring $R$, let $Q(R)$ be the total quotient ring of $R$
obtained from $R$ by inverting all non-zerodivisors.

\specialhead \noindent \boldLARGE 1. Grothendieck formula and zeta elements \endspecialhead

Let $X$ be a scheme of finite type over a finite field $\bF_q$. We assume $p$ is different from
$\text{char}(\bF_q)$.

In this $\S1$, we first review the formula (1.1.2) of
Grothendieck representing zeta functions of
$p$-adic sheaves on
$X$ by etale cohomology. We then show that those
zeta functions are recovered from $p$-adic zeta
elements (1.3.5).

{\bf 1.1. Zeta functions and etale cohomology groups in
positive characteristic case.} The Hasse zeta function
$\zeta(X, s) = \prod_{x \in |X|} \ (1 - \sharp
\kappa(x)^{-s})^{-1}$, where  $|X|$ denotes the set of all
closed points of $x$ and $\kappa(x)$ denotes the residue
field of
$x$, has the form
$\zeta(X, s) = \zeta(X/\bF_q, q^{-s})$ where
$$\zeta(X/\bF_q, u) = \prod_{x \in |X|} (1 -
u^{\text{deg}(x)})^{-1}, \ \ \ \ \ \text{deg}(x) = [\kappa(x) :
\bF_q].  \tag1.1.1$$

A part of Weil conjectures was that $\zeta(X/\bF_q, u)$ is a
rational function in $u$, and it was proved by Dwork and then
slighly later by Grothendieck. The proof of Grothendieck gives a
presentation of
$\zeta(X/\bF_q, u)$ by using etale cohomologyy.
More generally,
 for a finite extension $L$ of $\bQ_p$ and
for a constructible $L$-sheaf
$\cF$ on $X$, Grothendieck proved that the
L-function
 $L(X/\bF_q, \cF, u)$ has
the presentation
$$ L(X/\bF_q, {\cF}, u) = \prod_m \dt_L(1 - \varphi_q u \
; \  H^m_{et, c}(X \otimes_{\bF_q} {\bar \bF}_q ,
\cF))^{(-1)^{m-1}} \tag1.1.2$$
where $H_{et, c}^m$ is the etale cohomology with compact
supports and $\varphi_q$ is the action of the $q$-th power
morphism on $X$.

 In the case
$L =
\bQ_p = \cF$,
$\zeta(X/\bF_q, u) = L(X/\bF_q, \cF, u)$.
\medskip
{\bf 1.2. {\boldmath p}-adic zeta elements in positive characteristic case.} Determinants appear in the theory of
zeta functions as above, rather often. The regulator of a number field, which appears in the class number formula,
is a determinant.
 Such relation with determinant is well expressed
by the notion of ``determinant module".

If $R$ is a field, for an $R$-module $V$
of dimension $r$, $\dt_R(V)$ means the 1 dimensional
$R$-module $\wedge_{R}^r (V)$. For a bounded complex $C$ of
$R$-modules whose cohomologies $H^m(C)$ are finite
dimensional,
$\dt_{R}(C)$ means $\otimes_{m \in \bZ}
\ \{\dt_{R}(H^m(C))\}^{\otimes (-1)^m}$.

This definition is generalized to the definition of an
invertible
$R$-module $\dt_R(C)$ associated to a perfect complex
$C$ of
$R$-modules for a commutative ring $R$ (see [KM]).
 $\dt^{-1}_R(C)$ means the
inverse of the invertible module $\dt_R(C)$.

By a pro-$p$ ring, we mean a topological ring which
is an inverse limit of finite rings whose orders are
powers of $p$. Let
$\L$ be a commutative pro-$p$ ring. By a ctf $\L$-complex on
$X$, we mean a complex of $\L$-sheaves on
$X$ for the etale topology with constructible cohomology
sheaves and with perfect stalks. For a ctf $\L$-complex
$\cF$ on
$X$, $\RG_{et, c}(X, \cF)$ ($_c$ means with compact supports)
is  a perfect complex over $\L$.

For a commutative pro-$p$ ring $\L$ and for a ctf $\L$-complex
$\cF$ on
$X$, we define the
$p$-adic zeta element
$\zeta(X, \cF, \L)$ which is a $\L$-basis of
$\dt_{\L}^{-1}R\Gamma_{et, c}(X,
\cF)$. Consider the distinguished triangle
$$R\Gamma_{et, c}(X, \cF) \to R\Gamma_{et, c}(X
\otimes_{\bF_q}
\ {\bar
\bF}_q ,
\cF) @>1 - \varphi>> R\Gamma_{et, c}(X \otimes_{\bF_q} {\bar
\bF}_q ,
\cF). \tag1.2.1$$
Since det is multiplicative for distinguished triangles,
(1.2.1) induces an isomorphism
$$\dt_{\Lambda}^{-1} R\Gamma_{et, c}(X, \cF) \cong
\dt_{\L}^{-1} R\Gamma_{et, c}(X
\otimes_{\bF_q}  {\bar \bF}_q ,
\cF) \otimes_{\Lambda} \dt_{\L} R\Gamma_{et, c}(X
\otimes_{\bF_q}  {\bar \bF}_q ,
\cF) \cong \Lambda. \tag1.2.2$$
We define $\zeta(X, \cF, \L)$ to be the image of
$1
\in \Lambda$ in $\dt_{\L}^{-1} R\Gamma_{et, c}(X,
\cF)$ under (1.2.2). It is a $\L$-basis of the invertible
$\L$-module $\dt_{\L}^{-1} \RG_{et, c}(X,
\cF)$.

{\bf 1.3. Zeta functions and {\boldmath p}-adic zeta elements in positive characteristic case.} Let $L$ be a
finite extension of $\bQ_p$, let $O_L$ be the valuation ring of $L$, and let $\cF$ be a constructible $O_L$-sheaf
on $X$. We show that the zeta function $L(X/\bF_q, \cF_L, u)$ of the $L$-sheaf $\cF_L = \cF \otimes_{O_L} L$ is
recovered from a certain $p$-adic zeta element as in (1.3.5) below. Let
$$\L = O_L[[\Gal({\bar \bF}_q/\bF_q)]] = \varprojlim_{n}
O_L[\Gal(\bF_{q^n}/\bF_q)]. \tag1.3.1$$
 Let $s(\L)$ be the $\L$-module $\L$ which is regarded as a
sheaf on the etale site of $X$ via the natural action of
$\Gal({\bar \bF}_q/\bF_q)$. Then
$$H^m_{et, c}(X, \cF \otimes_{O_L} s(\L)) \cong
\varprojlim_{n} H^m_{et, c}(X
\otimes_{\bF_q} \bF_{q^n}, \cF) \tag1.3.2$$
where the transition maps of the inverse system are the trace
maps. From this, we can deduce that
$H^m_{et, c}(X,
\cF
\otimes_{O_L} s(\L))$ is a finitely generated $O_L$-module
for any $m$. Hence  we have $Q(\L) \otimes_{\L} \RG_{et, c}(X,
\cF
\otimes_{O_L} s(\L)) = 0$ and this gives an identificatition
canonical isomorphism
$$Q(\L) \otimes_{\L} \dt_{\L}^{-1} \RG_{et, c}(X, \cF
\otimes_{O_L} t(\L))
= Q(\L).  \tag1.3.3$$ Note
$$Q(\L) = Q(\varprojlim_n O_L[u]/(u^n - 1)) \supset Q(O_L[u]) =
L(u). \tag1.3.4$$
By a formal argument, we can prove the following (1.3.5)
(1.3.6) which show
\medskip
\ \  \ \  \ zeta function = zeta
element, \ \ \ \
 zeta value  = zeta element,
\medskip
\noindent respectively.
$$L(X/\bF_q, \cF_L, u) = \zeta(X, \cF \otimes_{O_L}
s(\L), \L) \ \ \text{in} \ \ Q(\L). \tag1.3.5$$
If $H^m_{et, c}(X, \cF_L) = 0$ for any $m$, $L(X/\bF_q, \cF_L, u)$
has no zero or pole at $u = 1$, and
$$L(X/\bF_q, \cF_L, 1) = \zeta(X, \cF, O_L) \ \ \text{in} \ \
L. \tag 1.3.6 $$

\specialhead \noindent \boldLARGE 2. Tamagawa number conjecture \endspecialhead

In 2.1, we describe the
generalized version of Tamagawa number conjecture. In 2.2
(resp. 2.3),  we consider $p$-adic zeta elements associated
to 1 (resp. 2) dimensional $p$-adic representations of
$\Gal({\bar \bQ}/\bQ)$, and their
relations to  (0.2) (resp. (0.3)).

{\bf
2.1. The conjecture.} Let
$X$ be a scheme of finite type over
$\1/p$. For a complex of sheaves $\cF$ on $X$ for the etale
topology, we define the compact support version
$R\Gamma_{et, c}(X,
\cF)$ of
$R\Gamma_{et}(X, \cF)$ as the mapping fiber of
$$R\Gamma_{et}(\1/p, Rf_!\cF) \to R\Gamma_{et}(\bR, Rf_!\cF)
\oplus
\RG_{et}(\bQ_p, Rf_!\cF).$$
where $f : X \to \Spec(\1/p)$.

It can be shown that for a commutative
pro-$p$ ring
$\L$ and for a ctf
$\L$-complex $\cF$ on $X$, $\RG_{et, c}(X,
\cF)$ is perfect.

The following is a generalized version of
the Tamagawa number conjecture [BK] (see [FP],  [Pe$_1$],
[Ka$_1$], [Ka$_2$]).
 In [BK], the
idea of Tamagawa number of motives was important, but it does
not appear explicitly in this version.

\proclaim{Conjecture}
To any triple $(X, \L, \cF)$ consisting of a scheme $X$ of
finite type over
$\1/p$, a commutative pro-$p$ ring $\L$, and a ctf
$\L$-complex on
$X$, we
can associate a
$\L$-basis
$\zeta(X, \cF, \L)$ of
$$\Delta(X, \cF, \L) =
\dt^{-1}_{\Lambda} R\Gamma_{et, c}(X, \cF),$$
which we call the $p$-adic zeta element associated to
$\cF$, satisfying the following conditions (2.1.1)-(2.1.5).

\noindent (2.1.1) If $X$ is a scheme over a finite field
$\bF_q$,
$\zeta(X, \cF, \L)$ coincides with the element defined in
$\S3.2$.

\noindent (2.1.2) (rough form) \  If $\cF$ is the $p$-adic
realization of a motive
$M$,
$\zeta(X, \cF, \L)$ recovers the complex value $\lim_{s
\to 0} s^{-e}L(M, s)$ where
$L(M, s)$ is the zeta function of $M$ and $e$ is the order of
$L(M, s)$ at $s = 0$.

\noindent (2.1.3) If $\L'$ is a pro-$p$ ring and $\L \to
\L'$ is a continuous homomorphism,
$\zeta(X, \cF \otimes_{\L}^L  \L', \L')$ coincides with the
image of
$\zeta(X,
\cF, \L)$ under
$\Delta(X, \cF
\otimes_{\L}^L \L', \L') \cong \Delta(X, \cF) \otimes_{\L}
\L'.$

\noindent (2.1.4)  For a distinguished triangle
$\cF' \to \cF \to \cF"$
with common $X$ and $\L$, we have
$$\zeta(X, \cF, \L) = \zeta(X, \cF', \L) \otimes
\zeta(X,
\cF", \L)
\
\
\text{in} \ \ \Delta(X, \cF, \L) = \Delta(X, \cF', \L)
\otimes_{\L} \Delta(X, \cF", \L).$$

\noindent (2.1.5)  If $Y$ is a scheme of finite type over $\1/p$ and $f : X \to Y$ is a separated morphism,
$$\zeta(Y, Rf_! \cF, \L) = \zeta(X, \cF, \L) \ \ \text{in} \
\
\Delta(Y, Rf_!
\cF, \L)
 = \Delta(X, \cF, \L).$$
\endproclaim

By this (2.1.5), the constructions of $p$-adic zeta elements
are reduced to the case $X = \Spec(\1/p)$.
How to formulate the part (4.1.2) of
this conjecture is reduced to the case of motives over
$\bQ$ by (2.1.5) and
 $L(M, s) =
L(Rf_!(M), s)$ (by philosophy of motives), where $f : X \to
\Spec(\1/p)$.

The conditions (2.1.3)-(2.1.5) are formal properties which
are analogous to formal properties of zeta functions. The
conditions (2.1.1) and (2.1.3)-(2.1.5) can be interpreted as

\noindent (2.1.6) \ The system $(X, \L, \cF) \mapsto \zeta(X, \cF,
\L)$ is an ``Euler system".

In fact,
let $L$ be a finite extension of $\bQ_p$, $S$ a finite set of
prime numbers containing $p$, and
let $T$ be a free $O_L$-module of finite rank
endowed with a continuous
$O_L$-linear action of $\Gal({\bar \bQ}/\bQ)$ which is
unramified outside $S$. For $m \geq 1$, let $R_m =
O_L[\Gal(\bQ(\zeta_m)/\bQ)]$ and let
\medskip
$z_m = \zeta_{R_m}(\1/p, j_{m,!}(T \otimes_{O_L}
s(R_m)), R_m) \in \dt_{R_m}^{-1}
\RG_{et, c}(\bZ[\zeta_m, \frac{1}{mS}], T).$

\noindent($j_m : \Spec(\bZ[\frac{1}{mS}]) \to \Spec(\1/p)) $.

Then the conditions (4.1.1) and
(4.1.3)-(4.1.5) tell that when $m$ varies, the
$p$-adic zeta elements $z_m$ form a system satisfying the
conditions of Euler systems formulated by Kolyvagin [Ko].

 We illustrate the relation (2.1.2) with zeta functions.

Let $M$ be a motive over $\bQ$, that is, a direct
summand of the motive $H^m(X)(r)$ for a proper smooth scheme
$X$ over $\bQ$ and for $r \in \bZ$, and
assume that $M$ is endowed with an action of a number field
$K$. Then the zeta function
$L(M, s)$ lives in
$\bC$, and the
$p$-adic zeta element lives in the world of
$p$-adic etale cohomology. Since these two worlds are too much
different in nature,  $L(M, s)$ and the $p$-adic zeta element
are not simply related.

However in the middle of
$\bC$ and the
$p$-adic world,

(a) there is a 1 dimensional $K$-vector space
$\Delta_K(M)$ constructed by the Betti
realization and the de Rham realization of $M$, and
$K$-groups (or motivic cohomology groups) associated to $M$.

Let $\infty$ be an Archimedean place of $K$. Then

(b) there is an
isomorphism
$$\Delta_K(M) \otimes_K K_{\infty} @>>\cong> K_\infty$$

\noindent constructed by Hodge theory and $K$-theory.

Let $w$ be a place of $K$ lying over $p$, let $M_w$  be the
representation of $\Gal({\bar \bQ}/\bQ)$ over $K_w$ associated to
$M$, and let $T$ be a $\Gal({\bar \bQ}/\bQ)$-stable
$O_{K_w}$-lattice in $M_w$. Then

(c) there is an isomorphism
$$ \align  &\Delta_K(M) \otimes_K K_w @ >>\cong>
 \dt^{-1}_{K_w} \RG_{et, c}(\1/p,
j_*M_w)\\
 =& \dt^{-1}_{O{K_w}} \RG_{et, c}(\1/p, j_*T) \otimes_{O_{K_w}}
K_w  \endalign $$

\noindent
where $j : \Spec(\bQ) \to \Spec(\1/p)$, constructed by
$p$-adic Hodge theory and $K$-theory.

See [FP] how to construct (a)-(c) (constructions require some
conjectures). The part (2.1.2) of the conjecture is:

(d) there exists a $K$-basis $\zeta(M)$ of
$\Delta_K(M)$ (called the rational zeta element associated to
$M$), which is sent to
$\lim_{s
\to 0} s^{-e}L(M, s)$ under the isomorphism  (b) where
$e$ is the order of $L(M, s)$ at $s = 0$, and to
$\zeta(\1/p, j_*T, O_{K_w})$ in  $ \dt_{K_{w}}^{-1} \RG_{et,
c}(\1/p, j_*M_w)$ under the isomorphism (c).

The existence of
$\zeta(M)$ having the relation with $\lim_{s
\to 0} s^{-e}L(M, s)$  was conjectured by Beilinson
[Be].

How zeta functions and $p$-adic zeta elements are related
is illustrated in the following diagram.
$$ \CD
\text{zeta functions side}  @.(\text{Betti})
@<\text{Hodge theory}<< (\text{de Rham}) @.\\
@.@A\text{regulator}AA     @AA\text{$p$-adic Hodge theory}A @.\\
@.(K\text{-theory}) @>>\text{Chern class}> (\text{etale}) @.
\text{$p$-adic zeta elements side.}
\endCD $$

We have the following
picture.
$$\CD
\text{automorphic rep} @<?<< \text{motives} @>{}>> p\text{-adic
Gal rep}
\\
  @V{}VV  @VV?V          @VV?V \\
\text{zeta functions} @. \text{rational zeta elements}
@. p\text{-adic zeta elememts}
\endCD
$$

The left upper arrow with a question mark shows the
conjecture that the map \{motives\} $\to$ \{zeta
functions\} factor through automorphic
representations, which is a
 subject
of non-abelian class field theory (Langlands
correspondences). As the other question marks indicate, we
do not know how to construct zeta elements in general, at
present.

{\bf 2.2. {\boldmath p}-adic zeta elements for 1 dimensional galois representations. } Let $\L$ be a commutative
pro-$p$ ring, and assume we are given a continuous homomorphism
$$\rho : \Gal({\bar \bQ}/\bQ) \to GL_n(\L) $$
which is unramified outside a finite set $S$ of prime
numbers $S$ containing $p$.
Let
$\cF = \L^{\oplus n}$ on which $\Gal({\bar \bQ}/\bQ)$ acts via
$\rho$, regarded as a sheaf on $\Spec(\bZ[\frac{1}{S}])$ for
the etale topology. We consider how to construct the $p$-adic
zeta element
$\zeta(\bZ[\frac{1}{S}], \cF, \L).$

In the case $n = 1$, we can use the ``universal objects" as
follows. Such $\rho$ comes from the canonical homomorphism

$$\rho_{\text{univ}} : \Gal({\bar \bQ}/\bQ) \to
GL_1(\L_{\text{univ}})
\ \ \text{where} \ \
\L_{\text{univ}} = \bZ_p[[\Gal(\bQ(\zeta_{Np^\infty})/\bQ)]]$$
for some $N \geq 1$ whose set of prime divisors coincide
with $S$ and for some continuous ring homomorphism
$\L_{\text{univ}} \to \L$.
We have $\cF \cong
\cF_{\text{univ}}
\otimes_{\L_{\text{univ}}} \L$. Hence
$\zeta(\bZ[\frac{1}{S}],
\cF,
\L)$ should be defined to be the image of
$\zeta(\bZ[\frac{1}{S}],
\cF_{\text{univ}}, \L_{\text{univ}})$. As is explained in
[Ka$_2$] Ch. I, 3.3,  $\zeta(\bZ[\frac{1}{S}],
\cF_{\text{univ}}, \L_{\text{univ}})$ is the pair of
the $p$-adic Riemann zeta function and a system of cyclotomic
units. Iwasawa main conjecture is regarded as the statemnet
that this pair is a $\L_{\text{univ}}$-basis of
$\Delta(\bZ[\frac{1}{S}],
\cF_{\text{univ}}, \L_{\text{univ}})$.

{\bf 2.3. {\boldmath p}-adic zeta elements for 2 dimensional Galois representations.} Now consider the case $n =
2$. The works of Hida, Wiles, and other people suggest that the universal objects $\L_{\text{univ}}$ and
$\cF_{\text{univ}}$ for 2 dimensional Galois representations in which the determinant of the action of the complex
conjugation is -1, are given by
$$\L_{\text{univ}} = \varprojlim_n \ p\text{-adic Hecke
algebras of weight 2 and of level} \ Np^n, $$
$$\cF_{univ} = \varprojlim_n \ H^1 \ \text{of modular curves
of level} \ Np^n.$$
Beilinson [Be] discovered ratinal
zeta elements in
$K_2$ of modular curves, and the images of these elements
in the etale cohomology under the Chern class maps become
$p$-adic zeta elements, and the inverse limit of these
$p$-adic zeta elements should be $\zeta(\bZ[\frac{1}{S}],
\cF_{\text{univ}}, \L_{\text{univ}})$ at least
conjecturally.
 By using this plan, the author obtained
$p$-adic zeta elements for motives associated to eigen cusp
forms of weight $\geq 2$, from Beilinson elements. Here it is
not yet proved that these $p$-adic zeta elements are
actually basis of $\Delta$, but it can be proved that they
have the desired relations with values $L(E, \chi,
1)$ and $L(f, \chi, r)$ ($1 \leq r \leq k-1$) for elliptic
curves over
$\bQ$ (which are modular by [Wi], [BCDT]) and for eigen cusp
forms of weight
$k
\geq 2$, and for Dirichlet charcaters $\chi$. Beilinson
elements are related in the Archimedean world to
$\lim_{s \to 0} s^{-1}L(E, \chi, s)$ for elliptic curves $E$
over
$\bQ$, but not related to $L(E,
\chi, 1)$. However since they become
universal (at least conjecturally) in the inverse limit in the
$p$-adic world, we can obtain from them $p$-adic zeta elements
related to $L(E, \chi, 1)$. Using these elements and
applying the method of Euler systems [Ko], [Pe$_2$], [Ru$_2$],
[Ka$_3$],  we can obtain the following results ([Ka$_4$]).

\proclaim{Theorem} Let $E$ be an elliptic curve over $\bQ$, let $N
\geq 1$, and let $\chi : \Gal(\bQ(\zeta_N)/\bQ) $\ \ $\cong
(\bZ/N\bZ)^{\times} \to \bC$ be a homomorphism. If $L(E, \chi, 1)
\neq 0$, the $\chi$-part of $E(\bQ(\zeta_N))$ and the $\chi$ part
of the Tate-shafarevich group of $E$ over $\bQ(\zeta_N)$ are
finite.
\endproclaim
The $p$-adic L-function $L_p(E)$ of $E$ is constructed from
the values
$L(E, \chi, 1)$.

\proclaim{Theorem} Let $E$ be an elliptic curve over $\bQ$
which is of good reduction at $p$.

(1) rank$(E(\bQ) \leq $ ord$_{s=1} L_p(E)$.

(2) Assume $E$ is ordinary at $p$. Let $\L =
\bZ_p[[\Gal(\bQ(\zeta_{p^{\infty}}/\bQ)]]$. Then the
$p$-primary Selmer group of
$E$ over
$\bQ(\zeta_{p^\infty})$ is $\L$-cotorsion and its
characteristic polynomial divides $p^nL_p(E)$ for some $n$.
\endproclaim

This result was proved by Rubin in the case of elliptic
curves with complex multiplication ([Ru$_1$]).

 As
described above, we can obtain
$p$-adic zeta elements of motives associated to eigen cusp
forms of weight
$\geq 2$. For such modular forms, we can prove the analogous
statement as the above (2).

Mazur and Greenberg conjectured that the
charcteristic polynomial of the above $p$-primary Selmer group
and  the $p$-adic L-function divide
each other.

\widestnumber\key{AAAAA}

\specialhead \noindent \boldLARGE References \endspecialhead

\ref \key Be \by Beilinson, A. \paper \rm Higher regulators and values of $L$-functions, {\it J. Soviet Math.}, 30
(1985), 2036--2070 \endref

\ref \key BK \by Bloch, S. and Kato, K. \paper \rm Tamagawa numbers of motives and L-functions, in The
Grothendieck Festschrift, 1, Progress in Math., 86, Burkhauser (1990), 333--400 \endref

\ref \key BCDT \by Breuil, C., Conrad, B, Diamond, F., Taylor, R. \paper \rm On the modularity of elliptic curves
over $\bold Q$: wild 3-adic exercises, {\it J. Amer. Math. Soc.}, 14 (2001), 834--939 \endref

\ref \key FP \by Fontaine, J. -M., and Perrin-Riou, B. \paper \rm Autour des conjectures de Bloch et Kato,
cohomologie Galoisienne et valeurs de fonctions L, Proc. Symp. Pure Math. 55, Amer. Math. Soc., (1994), 599--706
\endref

\ref \key Ka$_1$ \by Kato, K. \paper \rm Iwasawa theory and $p$-adic Hodge theory, {\it Kodai Math. J.}, 16
(1993), 1--31 \endref

\ref \key Ka$_2$ \by Kato, K. \paper \rm Lectures on the approach
to Iwasawa theory for Hasse-Weil $L$-functions via $B_{dR}$. I,
Arithmetic algebraic geometry (Trento, 1991), 50--163, Lecture
Notes in Math., 1553, Springer, Berlin (1993) \endref

\ref \key Ka$_3$ \by Kato, K. \paper \rm Euler systems, Iwasawa theory, and Selmer groups, {\it Kodai Math. J.},
22 (1999), 313--372 \endref

\ref \key Ka$_4$ \by Kato, K. \paper \rm $p$-adic Hodge theory and values of zeta functions of modular forms,
preprint \endref

\ref \key KM \by Knudsen, F., and Mumford, D. \paper \rm The projectivity of the moduli space of stable curves I,
{\it Math. Scand.}, 39, 1 (1976), 19--55 \endref

\ref \key Ko \by Kolyvagin, V. A. \paper \rm Euler systems, in The Grothendieck Festchrift, 2, Birkhouser (1990),
435--483 \endref

\ref \key Pe$_1$ \by Perrin-Riou, B. \paper \rm Fonction L
$p$-adiques des repr\'esentations p-adiques, \linebreak
Ast\'erisque 229 (1995) \endref

\ref \key Pe$_2$ \by Perrin-Riou, B. \paper \rm Systemes d'Euler p-adiques et th\'eorie d'Iwaswa, {\it Ann. Inst.
Fourier}, 48 (1998), 1231--1307 \endref

\ref \key Ru$_1$ \by Rubin, K. \paper \rm The ``main conjecture''
of Iwasawa theory for imaginary quadratic fields, {\it Inventiones
math.}, 103 (1991), 25--68 \endref

\ref \key Ru$_2$ \by Rubin, K. \paper \rm Euler systems, Hermann Weyl Lectures, {\it Annals of Math. Studies},
147, Princeton Univ. Press (2000) \endref

\ref \key Wi \by Wiles, A. \paper \rm Modular elliptic curves and Fermat's last theorem, {\it Ann. of Math.}, 141
(1995), 443--551 \endref

\enddocument